\documentclass[12pt, reqno]{amsart}
\usepackage{amsmath, amsthm, amscd, amsfonts, amssymb, graphicx, color}

\newtheorem{theorem}{Theorem}[section]
\newtheorem{lemma}[theorem]{Lemma}
\newtheorem{proposition}[theorem]{Proposition}

\theoremstyle{definition}

\newtheorem{example}[theorem]{Example}

\theoremstyle{remark}

\numberwithin{equation}{section}
\newcommand{\C}{\mathbb C}
\newcommand{\cala}{\mathfrak A}
\newcommand{\calb}{\mathfrak B}
\newcommand{\calh}{\mathfrak H}
\newcommand{\calm}{\mathfrak M}
\newcommand{\caln}{\mathfrak N}
\newcommand{\calx}{\mathfrak X}

\newcommand{\N}{\mathbb N}

\newcommand{\calk}{\mathfrak K}
\newcommand{\call}{\mathfrak L}

\newcommand{\ip}[1]{\langle#1\rangle}
\begin{document}
\title[A Kadison--Sakai type Theorem]{A Kadison--Sakai type Theorem}
\author[M. Mirzavaziri, M.S. Moslehian]{Madjid Mirzavaziri and Mohammad Sal Moslehian}
\address{Department of Mathematics, Ferdowsi University of Mashhad, P. O. Box 1159,
Mashhad 91775, Iran; \newline Banach Mathematical Research Group
(BMRG), Mashhad, Iran;
\newline Centre of Excellence in Analysis on
Algebraic Structures (CEAAS), Ferdowsi University of Mashhad, Iran.}
\email{mirzavaziri@math.um.ac.ir and madjid@mirzavaziri.com}
\email{moslehian@ferdowsi.um.ac.ir and
moslehian@ams.org}


\begin{abstract}
The celebrated Kadison--Sakai theorem states that every derivation
on a von Neumann algebra is inner. In this paper, we prove this
theorem for ultraweakly continuous $*$-$\sigma$-derivations, where
$\sigma$ is an ultraweakly continuous surjective $*$-linear mapping.
We decompose a $\sigma$-derivation into a sum of an inner
$\sigma$-derivation and a multiple of a homomorphism. The so-called
$*$-$(\sigma,\tau)$-derivations are also discussed.
\end{abstract}
\maketitle

\section{Introduction}

One of important questions in the theory of derivations is that
``When are all bounded derivations inner?'' Forty years ago,
R.V.~Kadison \cite{KAD} and S.~Sakai \cite{SAK} independently proved
that every derivation of a von Neumann algebra $\calm$ into itself
is inner. This was the starting point for the study of the so-called
bounded cohomology groups. This nice result can be restated as
saying that the first bounded cohomology group $H^1(\calm;\calm)$
(i.e. the vector space of derivations modulo the inner derivations)
vanishes. One may conjecture the vanishing of the higher order
bounded cohomology groups $H^n(\calm;\calm), n\geq 2$. It is
remarkable to note that if $\calm$ is a von Neumann algebra of type
I, ${\rm II}_1$, III, or is of type ${\rm II}_1$ and stable under
tensoring with the hyperfinite ${\rm II}_1$ factor, then
$H^n(\calm;\calm)=0$ and $H^n(\calm;{\mathcal B}(\calh))=0$ for all
$n\geq 1$. Not all type ${\rm II}_1$ factors are stable under
tensoring with the hyperfinite ${\rm II}_1$ factor, and so this
theorem leaves open the type ${\rm II}_1$ case. It is still unknown
whether $H^1(\calm,{\mathcal B}(\calh))$ is always zero; see
\cite{S-S1, S-S2}.

Suppose that $\cala$ and $\calb$ are two algebras, $\calx$ is a
$\calb$-bimodule and $\sigma:\cala\to \calb$ is a linear mapping, which is called the ground mapping. A
linear mapping $d:\cala\to \calx$ is called a $\sigma$-derivation if
$d(ab)=d(a)\sigma(b)+\sigma(a)d(b)$ for all $a,b\in \cala$. These
maps have been extensively investigated in pure algebra. Recently,
they have been treated in the Banach algebra theory (see \cite{B-V,
H-J-M-M, M-M1, M-M2, MOS, Z-T} and references therein).

There are some applications of $\sigma$-derivations to develop an
approach to deformations of Lie algebras which have many
applications in models of quantum phenomena and in analysis of
complex systems; cf. \cite{H-L-S}. A wide range of examples are as
follows:

(i) Every ordinary derivation of an algebra $\cala$ into an
$\cala$-bimodule $\calx$ is an $\iota_{\cala}$-derivation
(throughout the paper, $\iota_{\cala}$ denotes the identity map on
the algebra $\cala$);

(ii) Every endomorphism $\alpha$ on $\cala$ is a
$\frac{\alpha}{2}$-derivation;

(iii) For a given homomorphism $\rho$ on $\cala$ and a fixed
arbitrary element $X$ in an $\cala$-bimodule $\calx$, the linear
mapping $\delta(A)=X\rho(A)-\rho(A)X$ is a $\rho$-derivation of
$\cala$ into $\calx$ which is said to be an inner $\rho$-derivation.

(iv) Every point derivation $d:\cala\to\C$ at the character $\theta$
is a $(\theta,\theta)$-derivation.

In this paper, we investigate $\sigma$-derivations. We divide our
main work into three sections. The first section is devoted to prove
a Kadison--Sakai type theorem for $\rho$-derivations on von Neumann
algebras when $\rho$ is a homomorphism. In the next section, section
we briefly discuss an extension of one of our result to the
$(\sigma,\tau)$-derivations. In the last section, we decompose a
$\sigma$-derivation into a sum of an inner $\sigma$-derivation and a
multiple of a homomorphism. The importance of our approach is that
$\sigma$ is a linear mapping in general, not necessarily a
homomorphism.

A von Neumann algebra $\calm$ is an ultraweakly closed
$*$-subalgebra of $B(\calh)$ containing the identity operator $I$,
where $\calh$ is a Hilbert space. By the weak (operator) topology on
$B(\calh)$ we mean the topology generated by the semi-norms $T
\mapsto |\ip{T\xi,\eta}|\,\,(\xi, \eta \in \calh)$. We also use the
terminology ultraweak (operator) topology for the weak$^*$-topology
on $B(\calh)$ considered as the dual space of the nuclear operators
on $\calh$. When we speak of ultraweak continuity (weak continuity,
resp.) of a mapping between von Neumann algebras $\calm$ and $\caln$
we mean that we equipped both $\calm$ and $\caln$ with the ultraweak
topology (the weak topology, resp.). We refer the reader to
\cite{TAK} for undefined phrases and notations.

\section{Towards a Kadison--Sakai type theorem for $\sigma$-derivations}

We start our work with the following lemma.

\begin{lemma}[S. Sakai] \label{S} Let
$\rho: \calm \to \caln$ be an ultraweakly continuous $*$-epimorphism
and $\delta:\calm \to \caln$ be a $*$-$\rho$-derivation between von
Neumann algebras. Then there is a central projection $P\in\calm$ and
an $*$-isomorphism $\tilde{\rho}:\calm P\to\caln$ such that
$\delta:\calm P\to\caln$ is a $*$-$\tilde{\rho}$-derivation.
Moreover, $\delta(A)=0$ for each $A\in\calm(I-P)$. Also,
$\tilde{\rho}^{-1}\circ\delta$ is an ordinary derivation on $\calm
P$.
\end{lemma}
\begin{proof} Since $\rho$ is an ultraweakly
continuous $*$-homomorphism, its kernel is an ultraweakly closed
ideal of $\calm$. Hence there is a central projection $Q\in\calm$
such that $\ker\rho=\calm Q$. Set $P=I-Q$. Then
$\tilde{\rho}=\rho|{_{\calm P}}:\calm P \to \caln$ is a
$*$-isomorphism. We have
\[\delta(ABP)=\delta(APBP)=\delta(AP)\rho(AP)+\rho(AP)\delta(BP).\]
Hence $\delta$ is a $*$-$\tilde{\rho}$-derivation on $\calm P$.
Moreover, if $A=B(I-P)\in\calm(I-P)$ is a projection then
\[\delta(A)=\delta((BQ)^2)=\delta(BQ)\rho(BQ)+\rho(BQ)\delta(BQ)=0,\]
since $BQ\in\calm Q=\ker\rho$. The space $\calm(I-P)$ is a von
Neumann algebra, because it is $\ker\rho=\rho^{-1}(\{0\})$ and
$\rho$ is ultraweakly continuous. Thus $\calm(I-P)$ is generated by
its projection and so $\delta(A)=0$ for each $A\in\calm(I-P)$. The
last assertion is now obvious.
\end{proof}

\begin{theorem}\label{kad} If $\rho:\calm \to \caln$ is an ultraweakly continuous $*$-epimorphism and $\delta:\calm \to \caln$
be an ultraweakly continuous $*$-$\rho$-derivation between von
Neumann algebras, then there is an element $U\in\calm$ with
$\|U\|\leq\|\delta\|$ such that $\delta(A)=U\rho(A)-\rho(A)U$ for
each $A\in\calm$. In other words, $\delta$ is an inner
$\rho$-derivation.
\end{theorem}
\begin{proof}
The mapping $\tilde{\rho}^{-1}\circ\delta$ is an ordinary derivation
on $\calm P$, where $P$ is as in Lemma \ref{S}. Thus, by
\cite[Theorem 2.5.3]{SAK}, there is a $V\in\calm P$ with
$\|V\|\leq\|\tilde{\rho}^{-1}\circ\delta\|\leq\|\delta\|$ such that
$(\tilde{\rho}^{-1}\circ\delta)(A)=VA-AV$ for all $A\in\calm P$.
Thus
$\delta(A)=\tilde{\rho}(V)\tilde{\rho}(A)-\tilde{\rho}(A)\tilde{\rho}(V)$,
for all $A\in\calm P$. Putting $U=\tilde{\rho}(V)$ we have
$\delta(A)=U\rho(A)-\rho(A)U$ for all $A\in\calm P$. The later
equality is also valid for $A\in\calm(I-P)$, since both sides are 0
for these elements. Finally,
$\|U\|=\|\rho(V)\|\leq\|V\|\leq\|\delta\|$.
\end{proof}

The following example can illustrate the case that the ground
mapping is not a homomorphism.

\begin{example} Let $\calh$ be an infinite
dimensional Hilbert space with an orthonormal basis
$\{e_n\}_{n\in\N}$, and $\theta:{\mathcal B}(\calh)\to\C$ be a
non-zero character. Define $\rho:{\mathcal B}(\calh)\to{\mathcal
B}(\calh)$ by
\[\langle \rho(A)e_n,e_m\rangle=\left \{
\begin{array}{ll}
\theta(A) & \rm{if~ } n=m=1 \\
0&\rm{otherwise~}
\end{array}\right.\]
and $\delta:{\mathcal B}(\calh)\to{\mathcal B}(\calh)$ by
\[\langle \delta(A)e_n,e_m\rangle=\left \{
\begin{array}{ll}
\theta(A)&\rm{if~}(n,m)=(1,2)\rm{~or~}(n,m)=(2,1)\\
0&\rm{otherwise~}
\end{array}\right.\]
A simple verification by using the matrix form of $\rho(A)$ and
$\delta(A)$ shows that $\rho$ is a $*$-homomorphism which is not
surjective and $\delta$ is a $*$-$\rho$-derivation. Clearly,
$\delta(I)$ does not commute with $\rho(I)$.

Let there be a non-zero $*$-$\theta$-derivation $\zeta:{\mathcal
B}(\calh)\to\C$. Then $\delta:{\mathcal B}(\calh)\to{\mathcal
B}(\calh)$ defined by
\[\langle \delta(A)e_n,e_m\rangle=\left \{
\begin{array}{ll}
\zeta(A) & \rm{if~ } n=m=1 \\
\theta(A)&\rm{if~}(n,m)=(1,2)\rm{~or~}(n,m)=(2,1)\\
0&\rm{otherwise~}
\end{array}\right.\]
is a $*$-$\rho$-derivation which is not $\rho$-inner, is not a
$*$-homomorphism and does not decompose to a sum of a $\rho$-inner
derivation and a $*$-homomorphism.
\end{example}

\section{$(\sigma,\tau)$-derivations}

Assume that $\cala, \calb$ are $*$-algebras and
$\sigma,\tau:\cala\to \calb$ are $*$-linear mappings. By a
$*$-$(\sigma,\tau)$-derivation we mean a linear mapping $d:\cala\to
\calb$ preserving $*$ such that $d(ab)=d(a)\sigma(b)+\tau(a)d(b)$
for all $a,b\in \cala$. Obviously, a $*$-$\sigma$-derivation is a
$*$-$(\sigma,\sigma)$-derivation.

\begin{proposition} Let $\cala, \calb$ be $*$-algebras and $\sigma,\tau:\cala\to \calb$ be $*$-linear mappings.
Then any $*$-$(\sigma,\tau)$-derivation $d:\cala\to \calb$ is a
$*$-$\frac{\sigma+\tau}{2}$-derivation.\end{proposition}

\begin{proof} First we show that each $*$-$(\sigma,\tau)$-derivation is
a $*$-$(\tau,\sigma)$-derivation. We have
$d(ab)=d(b^*a^*)^*=(d(b)^*\sigma(a)^*+\tau(b)^*d(a)^*)^*=d(a)\tau(b)+\sigma(a)\tau(b)$.

Now we conclude that $d(ab)=\frac{1}{2}d(ab)+\frac{1}{2}d(ab)=
\frac{1}{2}(d(a)\sigma(b)+\tau(a)d(b))+\frac{1}{2}(d(a)\tau(b)+\sigma(a)d(b))=
d(a)\frac{\sigma+\tau}{2}(b)+\frac{\sigma+\tau}{2}(a)d(b)$.
\end{proof}

The previous theorem enables us to focus on $*$-$\sigma$-derivations
while we deal with $*$-algebras. In particular we obtain the
following generalization of Theorem \ref{kad}.

\begin{theorem} Let $\calm, \caln$ be von Neumann algebras, let
$\rho_1,\rho_2 :\calm \to \caln$ be ultraweakly continuous and norm
continuous $*$-homomorphisms such that $\rho_1+\rho_2$ be
surjective, and let $\delta:\calm\to\caln$ be an ultraweakly
continuous $*$-$(\rho_1,\rho_2)$-derivation. Then there is an
element $U$ in $\calm$ with
$\|U\|\leq\frac{\|\rho_1\|+\|\rho_2\|}{2}\|\delta\|$ such that
$\delta(A)=U\frac{\rho_1(A)+\rho_2(A)}{2}-\frac{\rho_1(A)+\rho_2(A)}{2}U$
for each $A\in{\calm}$.\end{theorem}

\section{A Kadison--Sakai type theorem for $\sigma$-derivations}

In this section, we aim to remove the condition of being a
homomorphism from the ground mapping and provide our main result on
the generalization of the Kadison--Sakai theorem. Indeed, we would
like to decompose $d$ as a direct sum of an inner
$*$-$\sigma$-derivation and a $*$-homomorphism.

Recall that if we have a Hilbert space direct sum
$\calh=\calk\oplus\call$, $P$ is the projection corresponding to
$\calk$ and $T\in {\mathcal B}(\calh)$, then the compression
$T_{\calk}$ to $\calk$ is the operator $T_{\calk}\in {\mathcal
B}(\calk)$ defined by $T_{\calk}(k)=PT(k)\,\,(k\in\calk)$.
Obviously, $PTP\leftrightarrow T_{\calk}$ is an isometric
$*$-isomorphism between von Neumann algebras $P{\mathcal B}(\calh)P$
and ${\mathcal B}(\calk)$. If $\sigma:\calm \to \calm$ is a linear
mapping such that $\sigma(A)P=P\sigma(A)$ for all $A\in \calm$, then
$\sigma_{\calk}:\calm \to P \calm P$, $\sigma_{\call}:\calm \to
(I-P)\calm (I-P)$ can be defined by $\sigma_{\calk}(A)=P\sigma(A)P,
\sigma_{\call}(A)=(I-P)\sigma(A)(I-P)$. Furthermore,
$\sigma=\sigma_{\calk}\oplus\sigma_{\call}$. Note that $P\calm P$
($(I-P)\calm(I-P)$, resp.) is a von Neumann subalgebra of $B(\calk)$
($B(\call)$, resp.).

\begin{example} Let $\calh$ be the separable Hilbert space with
orthonormal basis $\{e_n\}$. Define $\sigma:{\mathcal
B}({\calh})\to {\mathcal B}({\calh})$ by
\[\langle\sigma(A)e_n,e_m\rangle=\left \{
\begin{array}{ll}
\exp(\frac{\rm i}s-\frac{\rm i}r)\langle Ae_s,e_r\rangle & \rm{if~ } n=2s-1,m=2r-1 \\
\frac12\exp(\frac{\rm i}s-\frac{\rm i}r)\langle Ae_s,e_r\rangle
&\rm{if~ } n=2s,m=2r \\
0&\rm{otherwise~ }
\end{array}\right.\]
and $d:{\mathcal B}({\calh})\to {\mathcal B}({\calh})$ by
\[\langle d(A)e_n,e_m\rangle=\left \{
\begin{array}{ll}
(\frac1s-\frac1r)\exp(\frac{\rm i}s-\frac{\rm i}r)\langle Ae_s,e_r\rangle & \rm{if~ } n=2s-1,m=2r-1 \\
\exp(\frac{\rm i}s-\frac{\rm i}r)\langle Ae_s,e_r\rangle
&\rm{if~ } n=2s,m=2r \\
0&\rm{otherwise~ }
\end{array}\right.\]
Then $d$ is a $\sigma$-derivation. To show this let $E_{rs}$ be the
operator defined on $\calh$ by $$\langle E_{rs}e_n,e_m\rangle =
\delta_{rm}\delta_{sn},$$ where $\delta$ is the $\delta$-Kroneker
function. Then
\[\langle d(E_{rs})e_n,e_m\rangle=\left \{
\begin{array}{ll}
(\frac1s-\frac1r)\exp(\frac{\rm i}{s}-\frac{\rm i}{r}) & \rm{if~ } n=2s-1,m=2r-1\\
\exp(\frac{\rm i}s-\frac{\rm i}r)
&\rm{if~ } n=2s,m=2r \\
0&\rm{otherwise~ }
\end{array}\right.\]
and so
\[\langle d(E_{tu}E_{rs})e_n,e_m\rangle=\delta_{ur}\left \{
\begin{array}{ll}
(\frac1s-\frac1t)\exp(\frac{\rm i}s-\frac{\rm i}t) & \rm{if~ } n=2s-1,m=2t-1 \\
\exp(\frac{\rm i}s-\frac{\rm i}t)
&\rm{if~ } n=2s,m=2t \\
0&\rm{otherwise~ }
\end{array}\right.\]
On the other hand,
\[\langle \sigma(E_{rs})e_n,e_m\rangle=\left \{
\begin{array}{ll}
\exp(\frac{\rm i}s-\frac{\rm i}r) & \rm{if~ } n=2s-1,m=2r-1\\
\frac12\exp(\frac{\rm i}s-\frac{\rm i}r)
&\rm{if~ } n=2s,m=2r \\
0&\rm{otherwise~ }
\end{array}\right.\]
and so
\[\langle d(E_{tu})\sigma(E_{rs})e_n,e_m\rangle=\delta_{ur}\left \{
\begin{array}{ll}
(\frac1u-\frac1t)\exp(\frac{\rm i}s-\frac{\rm i}t) & \rm{if~ } n=2s-1,m=2t-1 \\
\frac12\exp(\frac{\rm i}s-\frac{\rm i}t)
&\rm{if~ } n=2s,m=2t \\
0&\rm{otherwise~ }
\end{array}\right.\]
and
\[\langle \sigma(E_{tu})d(E_{rs})e_n,e_m\rangle=\delta_{ur}\left \{
\begin{array}{ll}
(\frac1s-\frac1u)\exp(\frac{\rm i}s-\frac{\rm i}t) & \rm{if~ } n=2s-1,m=2t-1 \\
\frac12\exp(\frac{\rm i}s-\frac{\rm i}t)
&\rm{if~ } n=2s,m=2t \\
0&\rm{otherwise~ }
\end{array}\right.\]
This shows that
\[d(E_{tu}E_{rs})=d(E_{tu})\sigma(E_{rs})+\sigma(E_{tu})d(E_{rs})\]
and so $d$ is a $\sigma$-derivation.

\noindent Now let $\calk$ be the Hilbert space with orthonormal
basis $\{e_{2n-1}\}$ and $\call$ be the Hilbert space with
orthonormal basis $\{e_{2n}\}$. Then ${\calh}={\calk}\oplus{\call}$.
If $\sigma_{\calk}: {\mathcal B}(\calh)\to {\mathcal B}(\calk)$ and
$\sigma_{\call}: {\mathcal B}(\calh)\to {\mathcal B}(\call)$ are
defined as above, then $\sigma=\sigma_{\calk}\oplus\sigma_{\call}$.
By the same fashion we can write $d=d_{\calk}\oplus d_{\call}$. One
can easily verify that $d_{\calk}$ is a $\sigma_{\calk}$-derivation,
$d_{\call}$ is a $\sigma_{\call}$-derivation and $\sigma_{\calk}$ is
a $*-$homomorphism on unitary group, $\sigma_{\call}$ is $\frac12$
times a $*$-homomorphism, and
$d_{\calk}(A)=U\sigma_{\calk}(A)-\sigma_{\calk}(A)U$, where $U$ is
the operator defined by $\langle
Ue_{2n-1},e_{2m-1}\rangle=\frac{1}{n}\delta_{nm}$. Moreover
$d_{\call}=2\sigma_{\call}$. Note also that we have
$\sigma_{\calk}(I)=I_{\calk}$, $\sigma_{\call}(I)=\frac12
I_{\call}$, $d_{\calk}(I)=0$ and $d_{\call}(I)=I_{\call}$.
\end{example}

We will show that the situation described in the above example is
true in general.

\begin{proposition}\label{vw} Suppose that $\sigma: \calm\to\calm$ is an ultraweakly continuous surjective $*$-linear
mapping, $d:\calm\to\calm$ is an ultraweakly continuous
$*$-$\sigma$-derivation such that $d(I)$ is a central element of
$\calm$,
$\call_0=\bigcup_{B,C\in{\calm}}(\sigma(BC)-\sigma(B)\sigma(C))(\calh)$,
$\call$ is the closed linear span of $\call_0$ and
$\calk=\call^\perp$. Then

{\rm (i)} $
\calk=\bigcap_{B,C\in{\calm}}\ker(\sigma(BC)-\sigma(B)\sigma(C))$;

{\rm (ii)} If $P=P_{\calk}$ is the projection corresponding to
${\calk}$, then $\sigma(A)P=P\sigma(A)$ and $d(A)P=Pd(A)$ for all
$A\in\calm$;

{\rm (iii)} If we define $\delta:\calm\to P\calm P$ by
$\delta(A)=Pd(A)P$, $\rho:\calm\to P\calm P$ by
$\rho(A)=\sigma(A)P$, $\alpha:\calm\to (I-P)\calm (I-P)$ by
$\alpha(A)=(I-P)d(A)(I-P)$ and $\tau:\calm\to(I-P)\calm (I-P)$ by
$\tau(A)=(I-P)\sigma(A)(I-P)$, then $\delta$ is an ultraweakly
continuous $*$-$\rho$-derivation, $\alpha$ is an ultraweakly
continuous $*$-$\tau$-derivation, $d=\delta\oplus\alpha$ and
$\sigma=\rho\oplus\tau$. Moreover $\rho$ is an ultraweakly
continuous $*$-epimorphism;

{\rm (iv)} ${\calk}=\ker\delta(I)$ and
${\call}=\overline{\alpha(I)({\call})}$;

{\rm (v)} $\delta$ is an inner $\rho$-derivation.

{\rm (vi)} $\tau(I)=\frac12 I_{\call}$
\end{proposition}

\begin{proof} (i) For each $B,C\in{\calm}, h\in{\calh}$ and $k\in{\calk}$ we have
\begin{eqnarray*}
0&=&\langle (\sigma(BC)-\sigma(B)\sigma(C))(h),k\rangle
\\&=&\langle h,
(\sigma(BC)-\sigma(B)\sigma(C))^*(k)\rangle \\
&=&\langle h,(\sigma(C^*B^*)-\sigma(C^*)\sigma(B^*))(k)\rangle .
\end{eqnarray*}
Since $\calm$ is a $*$-subalgebra of ${\mathcal B}({\calh})$, we
infer that $(\sigma(BC)-\sigma(B)\sigma(C))(k)=0$ for each
$B,C\in{\calm}$ and $k\in{\calk}$. This shows that
${\calk}=\bigcap_{B,C\in{\calm}}\ker(\sigma(BC)-\sigma(B)\sigma(C))$.

(ii) Let $P=P_{\calk}$ be the projection corresponding to
${\calk}$. For each $B,C,A\in{\calm}$ and $k\in{\calk}$ we have
\begin{eqnarray*}
(\sigma(BC)-\sigma(B)\sigma(C))\sigma(A)(k)&=&(\sigma(BC)\sigma(A)-\sigma(B)
\sigma(C)\sigma(A))(k)\\
&=&(\sigma(BCA)-\sigma(BCA))(k)\\&=&0.
\end{eqnarray*}
This shows that $\sigma(A)({\calk})\subseteq{\calk}$ for each
$A\in\calm$. Since $\sigma$ is a $*$-linear mapping, we have
$\sigma(A)P=P\sigma(A)$ for all $A\in\calm$.

By using \cite[Lemma 2.2]{M-M2} we get
\[0=d(B)(\sigma(CA)-\sigma(C)\sigma(A))(k)=(\sigma(BC)-\sigma(B)\sigma(C))d(A)(k)\]
for all $k\in{\calk}$. This implies that $d(A)(k)\in{\calk}$ for all
$k\in{\calk}$. Hence $d(A)({\calk})\subseteq{\calk}$ for all
$A\in\calm$. Since $d$ is a $*$-linear mapping we conclude that
$d(A)P=Pd(A)$ for all $A\in\calm$.

(iii) Firstly, $\rho$ is a $*$-homomorphism. In fact, if
$B,C\in{\calm}$ then
\begin{eqnarray*}
\rho(BC)(k)&=&\sigma(BC)P(k)\\&=&\sigma(BC)(k)\\&=&\sigma(B)\sigma(C)(k)\\
&=&\sigma(B)\sigma(C)P^2(k)\\&=&\sigma(B)P\sigma(C)P(k)\\&=&\rho(B)\rho(C)(k)
\end{eqnarray*}
for all $k\in{\calk}$, and
\begin{eqnarray*}
\rho(BC)(\ell)=\sigma(BC)P(\ell)=0=\sigma(B)P\sigma(C)P(\ell)=\rho(B)\rho(C)(\ell).
\end{eqnarray*}
for all $\ell\in{\call}$. Hence $\rho$ is a homomorphism on
$\calm$. Moreover,
$\rho(A^*)=\sigma(A^*)P=\sigma(A)^*P=(P\sigma(A))^*=(\sigma(A)P)^*=\rho(A)^*,
(A\in{\calm})$.

Secondly, $\delta$ is a $\rho$-derivation, since if $A,B\in{\calm}$
then
\begin{eqnarray*}
\delta(AB)(k)&=&d(AB)P(k)\\&=&Pd(AB)(k)\\&=&Pd(A)\sigma(B)(k)+P\sigma(A)d(B)(k)\\
&=&P^2d(A)\sigma(B)P^2(k)+P^2\sigma(A)d(B)P^2(k)\\&=&Pd(A)PP\sigma(B)P(k)+P\sigma(A)PPd(B)P(k)\\
&=&(\delta(A)\rho(B)+\rho(A)d(B))(k),
\end{eqnarray*}
for all $k\in{\calk}$, and
\begin{eqnarray*}
\delta(AB)(\ell)&=&Pd(AB)P(\ell)\\&=&0\\&=&Pd(A)PP\sigma(B)P(\ell)+P\sigma(A)PPd(B)P(\ell)\\
&=&(\delta(A)\rho(B)+\rho(A)d(B))(\ell).
\end{eqnarray*}
for all $\ell\in{\call}$.

Similarly one can show that $\alpha$ is a $\tau$-derivation. It is
obvious that $d=\delta\oplus\alpha$ and $\sigma=\rho\oplus\tau$.

(iv) We have
\[ d(I)=d(I^2)=d(I)\sigma(I)+\sigma(I)d(I)=2\sigma(I)d(I)\,.\]
Since $\sigma$ is surjective, $d(I)$ is of the form $\sigma(E)$
for some $E\in \calm$. Now for each $k \in \calk$ we have
\[d(I)(k)=\sigma(E)(k)=\sigma(EI)(k)=\sigma(E)\sigma(I)(k)=d(I)\sigma(I)(k)\] and so
$d(I)(k)=d(I)(2\sigma(I)-I)(k)=0$ for each $k\in{\calk}$. Hence
$\delta(I)(k)=d(I)P(k)=0$ for each $k\in{\calk}$. Thus $k\subseteq
\ker\delta(I)$. On the other hand the compression operator
$\delta(I)$ belongs to ${\mathcal B}({\calk})$. Hence
$\ker\delta(I)\subseteq{\calk}$. Thus ${\calk}=\ker\delta(I)$.
Similarly one can show that ${\call}=\overline{\alpha(I)({\call})}$.

(v) By Theorem~\ref{kad} there is an element $U\in\calm$ such
that $\delta(A)=U\rho(A)-\rho(A)U$ for all $A\in\calm$. Hence
$\delta$ is an inner $\rho$-derivation.

(vi) Since $d(I)$ is a central element of $\calm$, one easily
deduce that $\alpha(I)=d(I)(I-P)$ is in the center of $\calm$.
Hence
\[\alpha(I)=\alpha(I^2)=\alpha(I)\tau(I)+\tau(I)\alpha(I)=2\tau(I)\alpha(I).\]
Thus $(2\tau(I)-I)\alpha(I)=0$. It follows from
$\overline{\alpha(I)({\call})}={\call}$ that $2\tau(I)=I_{\call}$.
\end{proof}

We can now establish a version of the Kadison--Sakai theorem as our
main theorem.

\begin{theorem}~\label{decom} Suppose that $\sigma:\calm\to\calm$ is
an ultraweakly continuous surjective $*$-linear mapping and
$d:\calm\to\calm$ is an ultraweakly continuous
$*$-$\sigma$-derivation such that $d(I)$ is a central element of
$\calm$. Then $\calh$ can be decomposed into
${\calk}\oplus{\call}$ and $d$ can be factored as the form
\[\delta\oplus 2Z\tau,\] where $\delta:{\calm}\to\calm$ is an
inner $*$-$\sigma_{\calk}$-derivation, $Z$ is a central element, and
$2\tau=2\sigma_{\call}$ is a $*$-homomorphism.
\end{theorem}

\begin{proof} It remains to show the result concerning $\call$. For $\ell\in{\call}$
we have
\[\alpha(A)(\ell)=\alpha(I)\tau(A)(\ell)+\tau(I)\alpha(A)(\ell)=
\alpha(I)\tau(A)(\ell)+\frac12\alpha(A)(\ell).\] Thus
\[\alpha(A)(\ell)=2\alpha(I)\tau(A)(\ell).\]
Putting $Z=\alpha(I)$ we get the result.
\end{proof}

\bibliographystyle{amsplain}

\end{document}